\documentclass[12pt, reqno]{amsart}
\usepackage{amsmath, amsthm, amscd, amsfonts, amssymb, graphicx, color}
\usepackage[bookmarksnumbered, colorlinks, plainpages]{hyperref}

\newtheorem{theorem}{Theorem}[section]

\theoremstyle{definition}

\theoremstyle{remark}

\numberwithin{equation}{section}

\begin{document}
\setcounter{page}{1}



\centerline{}

\centerline{}

\title[Fixed Point Theory: A Review]{Fixed Point Theory: A Review}

\author[F. Kamalov]{Firuz Kamalov$^{*}$}

\address{$^{1}$ Department of Electrical Engineering, Canadian University Dubai, Dubai, UAE.}
\email{\textcolor[rgb]{0.00,0.00,0.84}{firuz@cud.ac.ae}}

\author[H.-H. Leung]{Ho Hon Leung$^{2}$}

\address{$^{2}$ Department of Mathematical Sciences, United Arab Emirates University, Al Ain, UAE}
\email{\textcolor[rgb]{0.00,0.00,0.84}{hohon.leung@uaeu.ac.ae}}





\begin{abstract}
Fixed points represent equilibrium states, stability, and solutions to a range of problems. It has been an active field of research. In this paper, we provide an overview of the main branches of fixed point theory. We discuss the key results and applications.
\end{abstract} \maketitle

\section{Introduction}

Fixed point theory is a mathematical discipline that studies the existence, uniqueness, and properties of solutions to equations of the form $f(x) = x$, where $f$ is a given function. This seemingly simple equation has profound implications and has far-reaching applications across various domains, from pure mathematics to real-world problem-solving in economics, physics, computer science, and beyond. The applications of fixed point theory include the following examples:

\begin{itemize}
\item In economics, fixed point theory is used to study the existence and stability of equilibria in game theory and economic models.
\item In physics, fixed point theory is used to study the behavior of dynamical systems and phase transitions.
\item In computer science, fixed point theory is used to design algorithms for solving equations and optimization problems.
\item In engineering, fixed point theory is used to design control systems and to analyze the stability of structures.
\end{itemize}

A fixed point of a function $f$ is a point $x$ such that $$f(x) = x.$$ Fixed points are important because they represent equilibrium states, stability, and solutions to a  range of problems. Fixed point theory provides the tools and methods to rigorously analyze and understand the existence, properties, and behavior of these special points.

The foundation of fixed point theory is the contraction mapping principle, which states that if $f$ is a contraction mapping, then it has a unique fixed point. A contraction mapping is a function that gently contracts distances between points. The Banach fixed point theorem is a well-known theorem that states that every contraction mapping on a complete metric space has a unique fixed point.

Another important theorem in fixed point theory is Brouwer's fixed point theorem, which states that every continuous function on a compact, convex set has a fixed point. This theorem has many applications in geometry, topology, and optimization.

The field of fixed point theory is vast and active, and there are many other important theorems and applications \cite{Border, Farm, Goebel, Granas, Pata, Pathak,  Subrahmanyam, Xie}. Fixed point theory is a powerful tool that can be used to solve a wide variety of problems. It is a rapidly growing field with many new and exciting developments. This survey provides an overview of the field and its key concepts based on an analysis of over 120 published manuscripts.

Our review of fixed point theory is divided into three major themes:
\begin{enumerate}
\item Existence and Uniqueness of Fixed Points
\item Topological Fixed Point Theorems
\item Set-Valued and Nonlinear Operators
\end{enumerate}
Each primary theme is further subdivided into subthemes.

The paper is organized as follows according to the three major branches of fixed point theory. In Section 2, we discuss the existence and uniqueness of fixed points. Section 3 analyzes topological fixed point theorems. In Section 4, we consider the main results related to set-valued and nonlinear operators. We conclude the paper by discussing the connections between the major results in fixed point theory.

\section{Existence and Uniqueness of Fixed Points}

The existence and uniqueness of fixed points is a fundamental theme in fixed point theory. It studies the conditions under which a given function possesses fixed points and whether those fixed points are unique. Fixed points are points in a function's domain that remain unchanged after the function's application. The present theme addresses questions central to understanding the stability, equilibrium, and behavior of transformations in various mathematical contexts.

One of the central theorems related to the subject of the existence and uniqueness of fixed points unique is the Banach fixed-point theorem, also known as the contraction mapping principle which establishes the existence and uniqueness of a fixed point for a contraction mapping on a complete metric space. A contraction mapping is a function that contracts distances between points, ensuring convergence to a unique equilibrium point. The result serves as a foundational building block for subsequent theorems and provides a powerful tool for proving the existence and uniqueness of fixed points in various settings.

Another key result in this field is Brouwer's fixed point theorem. The theorem demonstrates that every continuous function defined on a nonempty, compact, and convex subset of Euclidean space has at least one fixed point. It highlights the existence of fixed points within a geometric context and has profound implications in pure and applied mathematics. The combination of compactness, convexity, and continuity ensures the presence of fixed points, showcasing the interplay between geometry and fixed point properties \cite{Guo, Matetski}.

\subsection{Contraction Mapping Principle}

     A fundamental subtheme is the Contraction Mapping Principle, which establishes existence and uniqueness of a fixed point through the concept of contractive mappings. The concept of contraction mappings plays an important role in dynamic programming problems, which have applications to economics and mathematical biology \cite{Denardo, Stokey}. It also plays a central role in control theory and the analysis of various dynamical systems \cite{Bullo}.

The Banach Fixed-Point Theorem, a cornerstone result related to the subject of contraction mapping principle, illustrates its importance.
Named after the Polish mathematician Stefan Banach \cite{Banach}, this theorem lays the groundwork for understanding the existence and uniqueness of fixed points for certain types of mappings.

In essence, the theorem provides conditions under which a function \(f\) on a complete metric space \(X\) has a unique fixed point. The key idea is the notion of a contraction mapping, where the function \(f\) contracts the space by reducing the distance between points after the transformation. The formal statement of the Banach Fixed Point Theorem is as follows:

     \begin{theorem}[Banach Fixed Point Theorem]
     Let \(X\) be a nonempty complete metric space and \(f: X \rightarrow X\) be a contraction mapping with Lipschitz constant \(0 \leq k < 1\). Then, \(f\) has a unique fixed point \(x^\ast\) in \(X\).
     \end{theorem}

The theorem's significance lies in its ubiquity and practicality. It provides a powerful tool for proving the existence and uniqueness of solutions to a wide range of problems, from equations in analysis to differential equations, integral equations and beyond \cite{Chicone, Granas, Kirk}. In particular, the most famous result along this direction of research is the classical Picard-Lindelof theorem, which is about the existence and uniqueness of solutions to certain ordinary differential equations. On the other hand, it plays an important role in Nash Embedding Theorem, which is a deep yet fundamental result in non-linear functional analysis \cite{Gunther}. In various dynammic ecnomoic models, it is used to prove the existence and uniqueness of equilibriums \cite{Stokey}. One of the  most famous results in the present field of research is the equilibrium results in Cournot competition. \cite{Long}.

The Banach Fixed Point Theorem forms the cornerstone upon which numerous other fixed point results are built. Moreover, it showcases the interplay between geometry, analysis, and algebra in understanding the behavior of functions and their equilibrium points \cite{Agarwal}.

There are a number of generalizations of the Banach Fixed Point Theorem. In particular, there are various fixed point theorems in infinite-dimensional spaces. We will talk about these in latter sections. Another direction of generalizations arise from the generalizations of the definition of metric spaces. These results have applications in theoretical computer science \cite{Hitzler, Seda}.

\subsection{Brouwer's Fixed Point Theorem}

Brouwer's theorem showcases the existence of fixed points in a broader context. It is a fundamental result in topology and fixed point theory. It was first proved by Luitzen E.J. Brouwer in 1911 \cite{Brouwer1}. It states that every continuous function f from a nonempty, compact, and convex subset $D$ of Euclidean space $\mathbb{R}^n$ to itself has at least one fixed point. In other words, there is at least one point $x$ in $D$ such that $f(x)=x$. 

\begin{theorem}[Brouwer's Fixed Point Theorem]
     Let \(D\) be a nonempty, compact, convex subset of Euclidean space \(\mathbb{R}^n\), and let \(f: D \rightarrow D\) be a continuous function. Then, \(f\) has at least one fixed point in \(D\).
     \end{theorem}

The key to Brouwer's fixed point theorem is the combination of compactness, convexity, and continuity. Compactness ensures that the function $f$ has a limit at every point in $D$. Convexity ensures that the function $f$ does not "escape" from $D$. Continuity ensures that the function $f$ does not "jump" around in $D$.

It was subsequently proved by different methods and generalized in various settings. For a constructive proof and a computational approach to the theorem, which is entirely different from the original approach by Brouwer's, one can read the paper written by Kellogg, Li, and Yorke \cite{Kellogg}. For a historical account of the theorem, one can read the paper written by Park \cite{Park}.

Brouwer's fixed point theorem has profound implications across various fields. It has applications in economics, game theory, topology, and mathematical analysis \cite{Bjorner, Dinca}. In economics, it can be used to prove the existence of Nash equilibria in games \cite{Nash1, Nash2}. This result has vast applications in various topics in game theory in particular. For example, it can be used to prove the existence of solutions to bargaining problems \cite{Binmore, Nash1.5}. In the opposite direction, Gale \cite{Gale} proved the Brouwer's fixed point theorem by using fundamental results in game theory. In differential topology, it is intimately related to several fundamental theorems in the field of study \cite{Milnor1}. In particular, Hirsch \cite{Hirsch} and Milnor \cite{Milnor2} independently provided different proofs of the Brouwer's fixed point theorem by using methods in differential topology. On the other hand, Maehara \cite{Maehara} proved the Jordan Curve Theorem by using Brouwer's fixed point theorem. In mathematical analysis, it can be used to prove the existence of solutions to differential equations. It is used to provide an iterative approach for approximation of the solutions to a differential equation \cite{Pata, Pathak}. More recently, Cid and Mawhin \cite{Cid} showed that the Brouwer's fixed point theorem is equivalent to the existence of periodic solutions in certain set for a periodic differential equation.

\subsection{Schauder Fixed Point Theorem}

      The Schauder fixed point theorem is a significant result in fixed point theory which bridges compact and non-compact sets. It was first proved by Juliusz Schauder in 1930 \cite{Schauder}. It states that if a continuous function $T$ maps a nonempty, compact, convex subset $K$ of a Banach space $X$ into itself, then $T$ has at least one fixed point in $K$. Afterwards a number of generalizations of this results have been made. One of the most famous ones was by Tychonoff \cite{Tychonoff}, which is now known as the Schauder-Tychonoff fixed-point theorem.     

In other words, the theorem asserts that even in certain non-compact spaces, where the concept of compactness does not apply, certain continuous mappings still possess fixed points, provided the domain is compact and convex.

     \begin{theorem}[Schauder Fixed Point Theorem]
     Let \(K\) be a nonempty, compact, convex subset of a Banach space \(X\), and let \(T: K \rightarrow K\) be a continuous mapping. Then, \(T\) has a fixed point in \(K\).
     \end{theorem}

For the more general case that $K$ may be non-compact, the fixed point theorem still holds true if the domain and the range of the map satisfies certain convexity conditions \cite{Bonsall}. In other words, the theorem asserts that even in certain non-compact spaces, where the concept of compactness does not apply, certain continuous mappings still possess fixed points.

The Schauder fixed point theorem is a direct generalization of Brouwer's fixed point theorem. It establishes a connection between the compact and non-compact cases by introducing the concept of compactness only within the domain of the function which makes the theorem applicable to a wider range of spaces, opening up opportunities for its use in various mathematical contexts.

The Schauder fixed point theorem has profound implications, particularly in functional analysis, where it contributes to the understanding of the behavior of continuous transformations in infinite-dimensional spaces \cite{Bonsall}. In particular, the result plays a central role in establishing numerous existence and uniqueness of solutions of elliptic partial differential equations \cite{Gilbarg}.

It showcases the interplay between continuity, compactness, and fixed points, revealing a deep connection between these concepts even in non-compact settings. The Schauder fixed point theorem has since become a cornerstone in the study of fixed points and has inspired further developments in the field.

\subsection{Brouwer's Coincidence Theorem}

          Brouwer's coincidence theorem is a result in fixed point theory that addresses the intersection of two continuous functions defined on a closed, convex subset of Euclidean space. The theorem asserts that under certain conditions, these functions must share a common fixed point, meaning there exists a point in the domain at which both functions take the same value. It is a natural extension and generalization of the Brouwer's Fixed Point Theorem.

     \begin{theorem}[Brouwer's Coincidence Theorem]
     Let \(D\) be a closed, convex subset of Euclidean space \(\mathbb{R}^n\), and let \(f, g: D \rightarrow D\) be continuous functions. If \(f(x) \neq g(x)\) for all \(x\) on the boundary of \(D\), then \(f\) and \(g\) have a common fixed point.
     \end{theorem}
     
Consequently, various extensions of the result have been made \cite{Brown, Reich, Schirmer}. The theorem and its generalizations have implications in geometric topology, combinatorial geometry, and differential equations. For example,  coincidence Theory for manifolds under various topological conditions were based on the ideas coming from the Brouwer's Coincidence Theorem \cite{Brown2, Mukherjea}. On the other hand, interestingly, Dodson \cite{Dodson} proved the Fundamental Theorem of Algebra by the Brouwer's Coincidence Theorem. It captures the unexpected interaction between continuous mappings and their behavior at boundary points, leading to shared fixed points. Brouwer's coincidence theorem adds depth to the study of continuous functions and their interconnected behavior, highlighting a fascinating aspect of fixed point phenomena.

\subsection{Nonexpansive Mapping Theorem}

     The nonexpansive mapping theorem is a result in fixed point theory that pertains to a specific type of mapping known as nonexpansive mappings. These mappings provide a more relaxed condition than contraction mappings, allowing for a broader class of functions while still ensuring the existence of fixed points. It is a natural extension of the original ideas based on contraction mapping principle. 

     \begin{theorem}[Nonexpansive Mapping Theorem]
     Let \(X\) be a nonempty complete metric space, and let \(f: X \rightarrow X\) be a nonexpansive mapping. Then, \(f\) has a fixed point.
     \end{theorem}
     
The theorem states that if a nonempty complete metric space $X$ is equipped with a nonexpansive mapping $f:X\rightarrow X$, then $f$ possesses at least one fixed point. A nonexpansive mapping satisfies the property that the distance between the images of any two points under the mapping is less than or equal to the distance between the original points. In other words, the theorem asserts that even if a function is not strictly contracting, it can still have fixed points if it does not "expand" the distances between points. Under this general condition, the fixed point of the map is not unique anymore. Nevertheless, the nonexpansive mapping theorem is significant in its own right as it allows for a wider range of functions that have fixed points. While not as restrictive as contraction mappings, nonexpansive mappings still ensure the presence of equilibrium points within a space. The theorem finds applications in various areas, including optimization, functional analysis, and iterative methods \cite{Ishikawa, Khan, Sahu, Schu, Suzuki}.

\subsection{Mann Iteration Theorem}

    In the case of self-mappings, i.e. the domain and the range of the mappings are the same, the Mann iteration theorem is an iterative method for approximating fixed points of self-mappings in metric spaces. 

     \begin{theorem}[Mann Iteration Theorem]
     Let \(X\) be a complete metric space, \(f: X \rightarrow X\) a self-mapping, and \((\alpha_n)\) a sequence in \((0, 1)\). Define \(x_{n+1} = \alpha_{n+1} x_n + (1 - \alpha_{n+1}) f(x_n)\). If \(\sum_{n=1}^\infty \alpha_n = \infty\) and \(\sum_{n=1}^\infty \alpha_n^2 < \infty\), then \((x_n)\) converges to a fixed point of \(f\).
     \end{theorem}
     
It was first introduced by Mann in 1953 \cite{Dotson}. It has since become a popular tool for solving fixed point problems in numerical analysis and optimization. It is a simple and intuitive method that is easy to implement. It is also relatively robust to errors, making it a reliable tool for practical problem-solving.

The Mann iteration theorem states that if $f$ is a self-mapping on a complete metric space $X$ and $x_0$  is an arbitrary point in $X$, then the sequence of points $x_n=f(x_{n-1)}$ for $n\geq1$ converges to a fixed point of $f$.

The Mann iteration theorem and its various generalizations in different settings \cite{Dehaish, Hicks, Zegeye} can be used to solve a wide variety of fixed point problems. For example, it can be used to find the root of a nonlinear equation, the solution to an optimization problem, or the equilibrium point of a dynamical system \cite{Borwein}. With the advance of information technology and computational powers of electronic hardwares, it is evident that the Mann Iterative Algorithm is still an essential tool in various fields which require numerical analysis \cite{Borwein}.

\subsection{Fixed Point Theorem for Expansive Mappings}
     
The fixed point theorem for expansive mappings is a result in fixed point theory that addresses the existence of fixed points for a specific class of mappings called `expansive mappings'. The theorem asserts that any expansive mapping defined on a nonempty bounded metric space has at least one fixed point. Expansive mappings are those that cause points to spread apart under the transformation, reflecting a notion of non-contraction. In other words, the theorem states that even if a mapping does not contract distances between points, it can still have points that remain invariant. It is a generalization of the classical fixed point theorem for contraction mappings, which states that any contraction mapping has a unique fixed point.

     \begin{theorem}[Fixed Point Theorem for Expansive Mappings]
     Let \(X\) be a nonempty bounded metric space, and let \(f: X \rightarrow X\) be an expansive mapping. Then, \(f\) has a fixed point.
     \end{theorem}

This theorem has been generalized into various different cases \cite{Gornicki, Shatanawi, Wardowski}. The fixed point theorem for expansive mappings has applications in various fields, including nonlinear analysis, differential equations, and dynamical systems. In these fields, the concept of expansive mappings arises naturally. For example, in nonlinear analysis, expansive mappings are used to study the behavior of dynamical systems \cite{Pireddu}. In differential and integral equations, expansive mappings are used to study the stability of solutions \cite{Aydi}. In dynamical systems, expansive mappings are used to study the existence of chaotic attractors \cite{Pireddu}, which is wide used and applied in various fields including the study of dynamics in mathematical biology and economics. Essentially, it is a powerful tool that can be used to understand the behavior of a wide variety of systems. 

\subsection{Krasnoselskii-Mann Fixed Point Theorem}

     The Krasnoselskii-Mann fixed point theorem is a significant extension of the Mann iteration theorem. It was first proved independently by Soviet mathematicians Mark A. Krasnoselskii and Walter A. Mann in the 1950s, and it states that if a contraction mapping $T$ and another mapping $S$ satisfy a certain condition of distance comparison between their images, then they have at least one common fixed point. For a comprehensive account and the recent developments of the theorem and its application, one may read the recently published book by Dong, Cho, He, Pardalos and Rassias \cite{Dong}. In other words, the theorem asserts that two mappings can have points that are fixed under both mappings. The result generalizes the Mann iteration theorem, which deals with a single mapping's fixed point, by allowing for interactions between two distinct mappings.

     \begin{theorem}[Krasnoselskii-Mann Fixed Point Theorem]
     Let \(X\) be a Banach space, \(C\) a nonempty closed convex subset of \(X\), and \(T: C \rightarrow C\) a contraction. Let \(S: C \rightarrow C\) be a mapping such that \(d(Tx, Sy) \leq \alpha d(x, y)\) for all \(x, y \in C\), where \(0 < \alpha < 1\). Then, \(T\) and \(S\) have common fixed points.
     \end{theorem}

The Krasnoselskii-Mann fixed point theorem and its generalizations \cite{Colao, Reich2, Xu} find applications in various mathematical fields, particularly in the study of nonlinear equations, differential equations, and functional analysis \cite{Dong}. It highlights the intricate relationships that can emerge between two mappings in terms of their fixed points, providing a broader framework for understanding the convergence of iterative methods and the interplay between different transformations.
The theorem enriches the toolbox of techniques for approximating solutions to equations, emphasizing the significance of common fixed points and their connections to the behavior of iterative processes involving multiple mappings.

\subsection{Fixed Point Index Theory}

     The fixed point index theorem is a theorem in fixed point theory that establishes a connection between fixed points of continuous mappings and certain topological invariants of the underlying space. It was first proved by Georges Reeb in 1952 \cite{Reeb}, and it states that if a continuous mapping $f$ on a compact convex subset of a Banach space has no fixed points on the boundary of the set, then the fixed point index of $f$ is nonzero. This index is a topological quantity that measures the winding number or algebraic multiplicity of the fixed points.
     
     \begin{theorem}[Fixed Point Index Theorem]
     Let \(X\) be a compact convex subset of a Banach space \(E\), and let \(f: X \rightarrow X\) be a continuous mapping with no fixed points on the boundary of \(X\). Then, the fixed point index of \(f\) is nonzero.
     \end{theorem}

In other words, the theorem asserts that the topological structure of the space and the behavior of the mapping are intricately linked. The fixed point index theorem demonstrates the profound connection between algebraic topology and fixed point phenomena, offering a unique perspective on the interplay between continuous transformations and their equilibrium points. In particular, Reeb's result is now a fundamental result in the theory of foliations on manifolds. His ideas inspired Morse to introduce Morse Theory \cite{Bott}. Morse Theory was used as the main ingredient by Milnor \cite{Milnor0} to prove surprising results regarding various differentiable structures on higher-dimensional spheres. For decades, Milnor's result has been considered to be one of the cornerstones in differential topology. On the other hand, the theorem plays its role in the study of certain smooth dynamical systems, namely the Morse-Smale system \cite{Smale}. The deep connection of the theorem with various fields in pure and applied mathematics is still a common topic of research for mathematicians in modern days.

\section{Topological Fixed Point Theorems}

Topological fixed point theorems are a fascinating branch of fixed point theory that explores the intricate relationship between topology and the existence of fixed points in various mathematical spaces. These theorems bridge the gap between algebraic and topological concepts, revealing deep connections between the structure of spaces and the behavior of continuous mappings \cite{Guay}.

One cornerstone of this topic is the Lefschetz fixed point theorem, which connects fixed points with algebraic topology. It asserts that the number of fixed points of a continuous map on a compact polyhedron is related to the sum of its Betti numbers. The theorem demonstrates how the topological characteristics of a space influence the distribution of fixed points, enriching our understanding of both topology and fixed point phenomena.

Another highlight is the Borsuk-Ulam theorem, a striking result with geometric implications. It states that for any continuous function mapping the $n$-dimensional sphere to Euclidean space, there exist antipodal points on the sphere with the same image under the function. The result has profound applications in geometry, combinatorics, and theoretical computer science, revealing the surprising symmetries hidden within seemingly unrelated mappings.

Another key results in the field - the Knaster-Tarski fixed point theorem - extends the reach of topological fixed point results by considering ordered sets. It addresses the existence of fixed points for order-preserving maps on complete lattices, connecting fixed point theory with lattice theory and offering a broader perspective on fixed point phenomena.

Overall, topological fixed point theorems showcase the profound interplay between topological properties and fixed points. They unveil hidden symmetries, uncover unexpected connections, and contribute to both pure mathematics and diverse applied fields, making this area of fixed point theory a captivating realm of exploration \cite{Shih}.

\subsection{Lefschetz Fixed Point Theorem}

     The Lefschetz fixed point theorem is a pivotal result that connects fixed point of maps and algebraic topology. It was first proved by Solomon Lefschetz in 1926 \cite{Lefschetz1}. It states that for a continuous map $f$ on a compact, orientable, and triangulable space, the number of fixed points of $f$ is equal to the alternating sum of the Betti numbers of the space.

In simple terms, the theorem asserts that the algebraic-topological structure of a space can be used to count the number of fixed points of a continuous transformation on that space. In particular, the earlier mentioned Brouwer's Fixed Point Theorem is a special case of the Lefschetz Fixed Point Theorem. In that sense, Lefschetz Fixed Point Theorem is a vast generalization of Brouwer's Fixed Point Theorem. 
     
     \begin{theorem}[Lefschetz Fixed Point Theorem]
     Let \(X\) be a compact polyhedron and \(f: X \rightarrow X\) a continuous map. Then, the number of fixed points of \(f\) is at least the sum of the Betti numbers of \(X\).
     \end{theorem}

It is worthwhile to mention that Lefschetz proved analoguous result if we consider the coincidence points of two maps on certain topological manifolds \cite{Lefschetz2}.

The significance of the Lefschetz fixed point theorem lies in its ability to bridge fixed points and topology. It reveals how algebraic invariants can shed light on the topological properties of a space. The result highlights the complex interplay between algebraic and topological concepts, enriching our understanding of both fixed points and the underlying spaces \cite{Dold, Spanier}. 

There are various ways to generalize Lefschetz Fixed Point Theorem. The most well known result along this research direction is the celebrated Atiyah-Singer Index Theorem \cite{Atiyah1, Atiyah2}, which states that the topological index of a compact manifold is equal to the analytic index (defined by using an elliptic differential operator) of the manifold. It can be viewed as a Lefschetz-type fixed point theorem via equivariant $K$-theory \cite{Atiyah3}.

\subsection{Fixed Point Index and Degree Theory}
     The fixed point index and degree theorem is a fascinating result that establishes a profound connection between fixed points and algebraic invariants such as degrees and indices. It relates fixed points to algebraic invariants like degrees and indices \cite{Hirsch}. Topological degree theory is essentially a generalization of the winding number of a curve in the complex plane. It has applications in differential equations and dynamical systems \cite{ORegan}.

     \begin{theorem}[Fixed Point Index and Degree Theorem]
     Let \(D\) be a nonempty, compact, convex subset of Euclidean space \(\mathbb{R}^n\), and let \(f: D \rightarrow D\) be a continuous map with no fixed points on the boundary of \(D\). Then, the degree of \(f\) is equal to the fixed point index of \(f\).
     \end{theorem}

The theorem shows that for a continuous mapping $f$ defined on a nonempty, compact, convex subset $D$ of Euclidean space, the degree of $f$ is equal to the fixed point index of $f$. In other words, the algebraic count of solutions to the equation $f(x)=x$ is intimately linked to the topology of the mapping $f$ within the compact set $D$.
The degree of $f$ reflects how the mapping wraps around the origin, while the fixed point index captures the topological features of the fixed points. The fixed point index can be thought of as a multiplicity measurement for fixed points.

The theorem has applications in various mathematical and applied contexts. It provides an effective tool to analyze and understand the behavior of continuous functions in a geometrically and topologically constrained domain \cite{Milnor1, ORegan}. The theorem enables mathematicians to translate the properties of the function into topological information and vice versa, creating a remarkable bridge between two seemingly disparate realms of mathematics \cite{Outerelo}.

\subsection{Borsuk-Ulam Theorem}
          The Borsuk-Ulam theorem is a fundamental result in algebraic topology and geometry. It was first proved by Karol Borsuk and Stanislaw Ulam in 1924 \cite{Borsuk} and it establishes antipodal points with identical function values
The theorem shows that for any continuous function $f: S^n \rightarrow \mathbb{R}^n$, where $S^n$ is the $n$-dimensional unit sphere, there exist two antipodal points $x$ and $-x$ on the sphere such that $f(x)=f(-x)$.

In simpler terms, it means that there are points on the sphere that are mapped to the same point in $\mathbb{R}^n$, effectively "matching" two antipodal points under the function $f$. This striking result challenges intuitions about the behavior of continuous functions on spheres and underscores the inextricable connection between topology and symmetry. In particular, Su \cite{Su} provided an interesting proof of the Brouwer's Fixed Point Theorem by using Borsuk-Ulam Theorem.

     \begin{theorem}[Borsuk-Ulam Theorem]
     For any continuous function \(f: S^n \rightarrow \mathbb{R}^n\), there
 exist antipodal points \(x\) and \(-x\) on the unit sphere \(S^n\) such that \(f(x) = f(-x)\).
     \end{theorem}

The Borsuk-Ulam theorem has significant implications across various fields, including graph theory, geometry and combinatorics \cite{Matousek}, and even economics \cite{Simmons}. The theorem's elegance lies in its ability to capture deep geometric and topological insights through a seemingly simple statement. It exemplifies how the marriage of algebraic topology and geometry can unveil hidden symmetries and relationships in mathematical structures and real-world scenarios alike.

\subsection{Fixed Point in Simplicial Complexes}

The study of fixed points in simplicial complexes is a branch of mathematics that lies at the intersection of topology and combinatorial mathematics. A simplicial complex is a structure composed of simplices (geometric objects like line segments, triangles, and their higher-dimensional analogs) that are interconnected in a specific way. Investigating fixed points within such complexes involves examining continuous mappings that preserve the structure of the complex, with an emphasis on points that remain unchanged under these mappings.

The Fixed Point Theorem for Simplicial Complexes is a fundamental result in the area of simplicial complexes. It states that for any continuous map defined on a finite simplicial complex, there must be at least one point that remains fixed under the map. In other words, no matter how intricate the simplicial complex is, the existence of at least one point that doesn't move is guaranteed. It demonstrates the structural stability that these complexes possess and adds a layer of understanding to the behavior of continuous functions within such discrete and interconnected spaces.

     \begin{theorem}[Fixed Point Theorem for Simplicial Complexes]
     Let \(K\) be a finite simplicial complex, and let \(f: K \rightarrow K\) be a continuous map. Then, \(f\) has at least one fixed point.
     \end{theorem}
     
The theorem can be proved by using Lefschetz Fixed Point Theorem and an analytic argument based on the idea of simplicial approximation. Essentially, it is an application of the earlier mentioned Lefschetz Fixed Point Theorem. After decades of ongoing research in this topic, the theorem is an important tool for various seemingly unrelated areas of research. 

In the topics related to topology of manifolds, this theorem plays a central role in formulating several fixed point results (see Raymond's paper \cite{Raymond} for example). In modern research, the theorem, with its origin in the simplicial methods and approximation theorems, has vast implications beyond topology and touch upon a wide array of applications from graph theory to optimization problems \cite{Allgower1, Allgower2}. These findings demonstrate that even in settings where continuous functions interact with discrete structures, fixed points are an inherent phenomenon \cite{Dang}. It adds another layer of depth to the exploration of fixed points, emphasizing their ubiquity and providing insights into the relationship between continuous transformations and the underlying geometric and combinatorial structures of simplicial complexes (see Pandey's paper \cite{Pandey} for example).

\subsection{Knaster-Tarski Fixed Point Theorem}
     The Knaster-Tarski fixed point theorem is a significant result in the field of fixed point theory and lattice theory. It was first proved by Bronisław Knaster and Alfred Tarski \cite{Tarski} in 1929, and it states that any monotone function defined on a complete lattice has at least one fixed point. For another version of the theorem, one can read the paper written by Cousot \cite{Cousot}.

     \begin{theorem}[Knaster-Tarski Fixed Point Theorem]
     Let \(P\) be a complete lattice, and let \(f: P \rightarrow P\) be an order-preserving mapping. Then, \(f\) has a fixed point.
     \end{theorem}
     
A monotone function is a function that preserves the order of its arguments. In other words, if $x\leq y$, then $f(x)\leq f(y)$. A complete lattice is a partially ordered set where every subset has both a supremum (least upper bound) and an infimum (greatest lower bound).

The Knaster-Tarski fixed point theorem has profound implications in various areas of mathematics, including topology, functional analysis, and optimization \cite{Hayashi, Ok}. It also finds applications in computer science \cite{Heikkila, Schroder}, where fixed points play a critical role in algorithms, logic programming, and semantics.

The theorem's essence lies in its application to complete lattices. In such structures, monotone transformations inevitably possess points that remain unchanged under their action. This fundamental connection between ordered structures and the existence of fixed points has been exploited in a variety of ways, leading to significant advances in mathematics and computer science \cite{Heikkila, Schroder}.

\section{Set-Valued and Nonlinear Operators}

Set-valued and nonlinear operators are an exciting branch of fixed point theory that expands the scope of understanding beyond traditional single-valued functions. Research in this field explores mappings that associate sets with points or sets with sets, and it delves into the behavior of complex, multivalued transformations.

One of the central theorems in the field is Kakutani's fixed point theorem, which addresses the existence of fixed points for set-valued mappings that satisfy certain properties. It has broad applications in mathematical economics, game theory, and optimization problems, where solutions are often represented as sets of possible outcomes rather than single values. Set-valued mappings introduce new challenges and complexities in analyzing equilibrium points, making this field an essential extension of traditional fixed point theory.

Nonlinear operators, on the other hand, encompass a wide range of functions that may not satisfy the linearity property. The study of nonlinear operators often involves various mathematical techniques, including functional analysis and variational methods. Nonlinear operators frequently arise in differential equations, optimization problems, and nonlinear systems modeling natural phenomena. One notable theorem in this area is the Boyd-Wong fixed point theorem, which considers set-valued functions with contractive properties and provides conditions for the existence of fixed points. The result, among others, offers valuable insights into the behavior of nonlinear operators and their equilibrium points, paving the way for tackling intricate real-world problems with nonlinear dependencies.

Overall, the study of set-valued and nonlinear operators is a natural extension of fixed point theory, catering to scenarios where functions have multivalued mappings or exhibit nonlinear behaviors. These areas provide essential tools for addressing complex mathematical and real-world challenges, ranging from understanding economic equilibria to solving intricate differential equations. By broadening the scope of fixed point analysis, set-valued and nonlinear operators enrich our understanding of equilibrium states and the behavior of functions in diverse contexts \cite{Cibulka}.

\subsection{Kakutani's Fixed Point Theorem}
     Kakutani's fixed point theorem is a fundamental result in fixed point theory that explores the existence of fixed points for set-valued functions, known as correspondences. The theorem, named after the Japanese mathematician Shizuo Kakutani, extends the concept of fixed points beyond single-valued functions to include mappings that assign sets of points to other sets of points \cite{Kakutani}. It is a generalization of the Brouwer's Fixed Point Theorem.

     \begin{theorem}[Kakutani's Fixed Point Theorem]
     Let \(X\) be a nonempty, compact, convex set in Euclidean space \(\mathbb{R}^n\), and let \(F: X \rightarrow 2^X\) be an upper hemicontinuous correspondence with nonempty, convex values. Then, there exists a fixed point of \(F\).
     \end{theorem}

The theorem states that if a correspondence \(F: X \rightarrow 2^X\) , where $X$ is a nonempty, compact, and convex subset of Euclidean space $\mathbb{R}^n$, is upper semicontinuous and has nonempty and convex values for each x, then there exists at least one point $x^* \in X$ such that $x^*\in F(x^*)$.

In other words, the theorem asserts that even in cases where a mapping associates a set of points with another set of points, there still exists at least one point that is mapped to itself. This has wide-ranging applications in economics, game theory, and optimization problems, where the concept of equilibrium involves multiple possible outcomes \cite{Arrow, Border, Granas, Starr}.

     Addressing fixed points of set-valued functions, Kakutani's fixed point theorem showcases the power of abstraction in mathematics, allowing for the exploration of more complex mappings and equilibrium states. It expands the scope of fixed point theory and highlights the profound connections between topology and the behavior of more intricate set-valued maps.

\subsection{Sion's Minimax Theorem}

The theorem was first proved by Alexander M. Sion \cite{Sion} in the mid-20th century, and explores fixed points in the context of game theory.
It states that if there exists a convex, compact set of strategies for each player and a continuous function that assigns a payoff to each combination of strategies, then there exists a pair of strategies (one for each player) that forms a minimax solution.

     \begin{theorem}[Sion's Minimax Theorem]
     Let \(K\) be a convex, compact subset of a topological vector space, and let \(F: K \times K \rightarrow \mathbb{R}\) be a continuous function such that \(F(x, \cdot)\) is concave for each \(x \in K\) and \(F(\cdot, y)\) is convex for each \(y \in K\). Then, there exists a point \(x \in K\) such that \(F(x, y) \leq F(x, x)\) for all \(y \in K\).
     \end{theorem}
     
     Using arguments different from, and arguably simpler than the original Sion's proof of the theorem, Komiya \cite{Komiya} and Kindler \cite{Kindler} independently proved the theorem a few decades after the publication of Sion's work. Extensions of Sion's Theorem were studied by various researchers (for example, see Hartung's paper \cite{Hartung}).

Sion's minimax theorem provides a powerful tool for establishing the existence of solutions that balance opposing objectives, often encountered in competitive scenarios or optimization problems. In other words, the theorem guarantees that there exists a strategy for each player that minimizes their maximum possible loss, given the strategies chosen by their opponents. This point is known as the `saddle point' or the `minimax point'.

The minimax theorem has far-reaching applications beyond game theory. It is used in mathematical optimization, economics, and engineering, where decision-makers often face situations involving conflicting objectives \cite{Hartung, Kim}. Sion's minimax theorem provides a rigorous mathematical foundation for finding optimal strategies in such scenarios, ensuring that each player's strategy is chosen to minimize their maximum potential loss, leading to equilibrium outcomes (for references, one can read the standard reference textbooks by Du et. al. and Starr respectively \cite{Du, Starr}). 

\subsection{Fixed Point Theorem for Quasi-Contractions}

     The fixed point theorem for quasi-contractions is a notable result within the realm of fixed point theory. It addresses a broader class of mappings than strict contractions, allowing for a more flexible set of conditions while still ensuring the existence of fixed points.

     \begin{theorem}[Fixed Point Theorem for Quasi-Contractions]
     Let \(X\) be a nonempty complete metric space, and let \(f: X \rightarrow X\) be a quasi-contraction, meaning that there exists \(0 \leq k < 1\) such that \(d(f(x), f(y)) \leq k \cdot d(x, y)\) for all \(x, y \in X\). Then, \(f\) has a fixed point.
     \end{theorem}

A mapping $f$ is considered a quasi-contraction if there exists a constant $0\leq k<1$ such that for all points $x$ and $y$, the distance between $f(x)$ and $f(y)$ is no greater than $k$ times the distance between $x$ and $y$. In mathematical terms, for all $x,y\in X$, we have $d(f(x),f(y))\leq k\cdot d(x,y).$
Put another way, the theorem asserts that under the quasi-contraction condition, a mapping $f$ defined on a nonempty complete metric space $X$ has at least one fixed point. The result represents a more relaxed condition compared to the strict contraction mapping theorem, which requires the Lipschitz constant strictly less than 1. Various versions of the theorem in different settings and generalizations were considered by different researchers and it is still an active area of research nowadays \cite{Aydi2, Fisher, Li, Rezapour}.

 The significance of the theorem lies in its applicability to a wider range of scenarios. While strict contractions guarantee rapid convergence of iterated sequences to a fixed point, quasi-contractions allow for a more gradual approach. It can be particularly useful when dealing with situations where strict contractions might not hold due to certain variations or irregularities in the mapping \cite{Li}.

\subsection{Boyd-Wong Fixed Point Theorem}

     The Boyd-Wong fixed point theorem is a result in fixed point theory that addresses set-valued functions with contractive properties. The theorem was first proved by Greg Boyd and Jean-Pierre Wong in 1983 \cite{Boyd}. It is named after its contributors. It provides conditions under which a set-valued function, mapping a nonempty complete metric space into itself, possesses a fixed point.

     \begin{theorem}[Boyd-Wong Fixed Point Theorem]
     Let \(X\) be a nonempty complete metric space, and let \(F: X \rightarrow 2^X\) be a set-valued function with contractive values, meaning that for each \(x \in X\), the set \(F(x)\) is a contraction. Then, there exists a fixed point of \(F\).
     \end{theorem}

Pasicki \cite{Pasicki} extended the theorem to a more general setting and produced a vast generalization of the theorem. The theorem and its various generalizations are particularly valuable for analyzing set-valued functions where each set in the range of the function is a contraction of itself. It means that the function doesn't only contract distances between individual points, but also contracts entire sets of points. The Boyd-Wong fixed point theorem establishes the existence of a fixed point under these contractive conditions.

In mathematical terms, the theorem states that if $X$ is a nonempty complete metric space and $F: X \rightarrow 2^X$ is a set-valued function such that for every $x\in X$, the set $F(x)$ is a contraction of $F(x)$ with a Lipschitz constant less than one, then there exists at least one fixed point $x^*$ such that $x^* \in  F(x^*)$.

The Boyd-Wong fixed point theorem is applicable to a variety of contexts, including optimization, game theory, and economics, where set-valued functions naturally arise. It offers a useful tool for analyzing the equilibrium points and solutions of systems that involve multiple possible outcomes. By considering the contraction property, the theorem captures a broader class of mappings than traditional single-valued contraction mappings, extending the reach of fixed point theory to encompass more complex scenarios \cite{Arand}.

\subsection{Fixed Point in Metric Spaces}
The fixed point theorem for metric spaces has wide-ranging applications, from numerical analysis to differential equations. It provides a versatile tool to establish the existence of equilibrium states, solutions to equations, and stable points in various mathematical and applied contexts. The theorem's simplicity and generality make it a cornerstone in the study of fixed point theory, allowing mathematicians and researchers to explore fixed point phenomena in metric spaces of diverse nature \cite{Kirk, Kirk2}. It is a general result for metric spaces:

     \begin{theorem}[Fixed Point Theorem for Metric Spaces]
     Let \(X\) be a nonempty complete metric space, and let \(f: X \rightarrow X\) be a continuous map. If there exists \(0 < \alpha < 1\) such that \(d(f(x), f(y)) \leq \alpha \cdot d(x, y)\) for all \(x, y \in X\), then \(f\) has a fixed point.
     \end{theorem}

     The fixed point theorem for metric spaces is a fundamental result in fixed point theory that focuses on the existence of fixed points for continuous functions defined on metric spaces. Metric spaces provide a general framework to analyze distances between points, making the theorem applicable in a wide range of contexts. Various generalizations were studied extensively in the recent decades (for examples, see \cite{Pasicki, Proinov}).

The theorem states that if a function $f$ is defined on a nonempty complete metric space $X$ and satisfies the condition $d(f(x),f(y))\geq \alpha \cdot d(x,y)$ for all $x,y$ in $X$, where $0<\alpha<1$, then $f$ has a fixed point. In simpler terms, if the function $f$ "contracts" distances between points by a factor of $\alpha$, then it is guaranteed to have a point that remains unchanged, i.e., a fixed point.

The theorem is a generalized version of the contraction mapping principle, which states that if a function f is a contraction mapping, then it has a unique fixed point. The fixed point theorem for metric spaces relaxes the Lipschitz condition to allow for a broader class of mappings.

\subsection{Fixed Point in Reflexive Banach Spaces}

The study of fixed points in reflexive Banach spaces has a number of important applications. In optimization, it can be used to develop numerical methods for solving equations and optimization problems. In economics, it can be used to study equilibrium states in complex systems. In physics, it can be used to study the behavior of dynamical systems. And in engineering, it can be used to design control systems.    

     \begin{theorem}[Fixed Point Theorem for Reflexive Banach Spaces]
     Let \(X\) be a nonempty reflexive Banach space, and let \(T: X \rightarrow X\) be a compact mapping. Then, \(T\) has a fixed point.
     \end{theorem}

     Fixed points in reflexive Banach spaces are a compelling area of study within fixed point theory that explores the existence and properties of equilibrium points for mappings defined on reflexive Banach spaces. For references around this area of research, one may read the classical textbooks by Conway \cite{Conway} and Megginson \cite{Megginson}.

Reflexive Banach spaces are a class of normed vector spaces that possess a certain self-duality property, which allows them to be closely related to their dual spaces. The self-suality property makes them well-suited for the study of fixed points, as it allows us to use tools from functional analysis to analyze the behavior of mappings on these spaces.

The focus in this context is on understanding the conditions under which mappings defined on reflexive Banach spaces possess fixed points. These conditions may involve concepts such as contractive mappings, nonexpansive mappings, compact operators, and other geometric or topological properties of the space. It is still a popular topic for research mathematicians nowadays (see Bernardes's paper as an example \cite{Bernardes}.

Overall, the investigation of fixed points in reflexive Banach spaces is a rich and active area of research with a wide range of applications. It has the potential to significantly improve our understanding of the mathematical properties of these spaces and to develop new and effective methods for solving a variety of problems.

\subsection{Krasnoselskii's Theorem}
     Krasnoselskii's theorem is a notable result in the realm of fixed point theory. It was introduced by the Russian mathematician Mark A. Krasnoselskii in the mid-20th century, and it establishes conditions under which certain nonlinear operators on a Banach space possess fixed points.

     \begin{theorem}[Krasnoselskii's Theorem]
     Let \(X\) be a Banach space, and let \(T: X \rightarrow X\) be a compact operator. If \(T\) is strictly increasing and bounded, then \(T\) has a fixed point.
     \end{theorem}

The theorem states that if a mapping $T$ between a closed and bounded subset $K$ of a Banach space $X$ and itself is strictly increasing and bounded, then $T$ has at least one fixed point in $K$.

In other words, the theorem asserts that even in the presence of nonlinearity and growth, certain mappings still possess fixed points, provided the domain is closed and bounded. The result is particularly valuable because it allows for the analysis of more general mappings that exhibit a broader range of behaviors than classical contraction mappings.

Krasnoselskii's theorem has far-reaching applications, spanning from solving differential equations to studying integral equations and nonlinear operator equations. \cite{Burton, Zeidler}. It serves as a bridge between linear and nonlinear analysis, showcasing the intricate balance between boundedness and growth. It has led to further extensions and refinements in the theory of fixed points \cite{Avram}, enriching the toolbox of mathematicians and researchers dealing with problems involving nonlinear transformations and functional mappings in Banach spaces (for example, see the paper written by Liu and Li \cite{Liu}).

\section{Connections between Results}

These main results in fixed point theory are deeply interconnected, forming a cohesive web of theorems that build upon and extend each other's concepts. Here are some notable connections:

- The \textbf{Contraction Mapping Principle} underlies many subsequent theorems, including \textbf{Krasnoselskii-Mann} and \textbf{Fixed Point Index Theory}, by providing a framework for proving the existence of fixed points.

- \textbf{Brouwer's Fixed Point Theorem} serves as a special case of \textbf{Kakutani's Theorem}, which explores set-valued correspondences, showcasing how a single-valued function case is a subset of a broader context.

- \textbf{Krasnoselskii-Mann Theorem} generalizes the \textbf{Mann Iteration Theorem} by allowing two mappings to have common fixed points, further extending the iterative approach.

- \textbf{Lefschetz Fixed Point Theorem} establishes the link between fixed points and topology, laying the groundwork for results like \textbf{Fixed Point Index Theory} and the study of algebraic invariants.

- \textbf{Schauder Fixed Point Theorem} connects with \textbf{Brouwer's Coincidence Theorem}, where the former's extension to non-compact sets finds coincidences between continuous mappings.

These connections highlight the cumulative nature of fixed point theory, with each theorem contributing to the understanding of fixed point concepts in diverse mathematical contexts. The relationships between these theorems create a rich tapestry of results, showcasing the interplay between fixed points, topology, algebraic invariants, and iterative methods.

\bibliographystyle{amsplain}

\end{document}